\documentclass{article}
\usepackage{amsfonts,amsmath,amsthm,pgf}

\newcommand\ad{\operatorname{ad}}

\newcommand\CC{\mathbb{C}}
\newcommand\RR{\mathbb{R}}
\newcommand\dexp{\operatorname{dexp}}
\newcommand\Id{\mathrm{I}}

\newtheorem{theorem}{Theorem}
\newtheorem{lemma}[theorem]{Lemma}
\newtheorem{conjecture}[theorem]{Conjecture}
\theoremstyle{definition}
\newtheorem{example}{Example}

\begin{document}

\title{Convergence of the Magnus series}
\author{Per Christian Moan\footnote{Centre of Mathematics for
  Applications, University of Oslo, PO Box 1053 Blinders, \hbox{NO-0316}
  Oslo, Norway.} 
  \and Jitse Niesen\footnote{Department of Mathematics, La Trobe
  University, Melbourne, Victoria 3086, Australia. email:
  j.niesen@latrobe.edu.au}}
\maketitle

\begin{abstract}
The Magnus series is an infinite series which arises in the study of
linear ordinary differential equations. If the series converges, then
the matrix exponential of the sum equals the fundamental solution of
the differential equation. The question considered in this paper is:
When does the series converge? The main result establishes a sufficient
condition for convergence, which improves on several earlier results.
\end{abstract}

\section{Introduction}

The Magnus series is an infinite series which arises in the study of
linear ordinary differential equations of the form $y' = A(t)\,y$,
$y(0) = y_0$, where $y(t)$ is a vector and $A(t)$ is a matrix. We
assume throughout this paper that $A(t)$ is a real-valued matrix, even
though the Magnus series is valid in more general settings.

The fundamental solution is defined by
\begin{equation}
\label{ode}
Y' = A(t) \, Y, \qquad Y(0) = \Id.
\end{equation}
If $A(t)$ is constant, then the solution of~\eqref{ode} is
given by the matrix exponential $Y(t) = e^{At}$. This suggests the
ansatz $Y(t) = e^{\Omega(t)}$, where $\Omega$ is a matrix function to
be determined, for the nonautonomous equation~\eqref{ode}. It turns out
that $\Omega(t)$ satisfies the differential equation
\begin{equation}
\label{omega-eqn}
\Omega' = \dexp^{-1}_{\Omega} \big(A(t)\big) 
= \sum_{k=0}^\infty \frac{B_k}{k!} \ad_\Omega^k \big(A(t)\big),
\end{equation}
where $B_k$ denote the Bernoulli numbers ($B_0=1$, $B_1=-\frac12$,
$B_2=\frac16$, $B_3=0$, etc.) and $\ad_\Omega$ is the adjoint operator,
defined by
$$
\ad_\Omega(A) = [\Omega,A] = \Omega A - A \Omega.
$$
Magnus~\cite{magnus54oes} applied Picard iteration on~\eqref{omega-eqn}
to find an infinite series for~$\Omega$:
\begin{align}
\Omega(t) = \int_0^t A(\tau) \,d\tau 
&- \frac12 \int_0^t \int_0^{\tau_1} \big[ A(\tau_2), A(\tau_1) \big]
d\tau_2 \, d\tau_1 \notag \\ 
&+ \frac14 \int_0^t \int_0^{\tau_1} \int_0^{\tau_2} \Big[ \big[ 
A(\tau_3), A(\tau_2) \big], A(\tau_1) \big] d\tau_3 \,d\tau_2
\,d\tau_1 \label{MagnExp} \\ 
&+ \frac1{12} \int_0^t \int_0^{\tau_1} \int_0^{\tau_1} \Big[
A(\tau_3), \big[ A(\tau_2), A(\tau_1) \big] \Big] d\tau_3 \,d\tau_2
\,d\tau_1 + \cdots . \notag
\end{align}
This series has since come to be called the \emph{Magnus series}. We
can write it as
$$
\Omega(t) = \Omega_1(t) + \Omega_2(t) + \Omega_3(t) + \Omega_4(t) +
\cdots, 
$$
where the term $\Omega_n(t)$ is a sum of $n$-fold integrals of $n-1$
nested commutators. Explicit expressions for the $\Omega_n(t)$ are
given by Bialynicki-Birula, Mielnik and
Plaba\'nski~\cite{bialynicki-birula69eso}, Chacon and
Fomenko~\cite{chacon91rff}, and Iserles and
N\o{}rsett~\cite{iserles99oso}.

The Magnus series can be used to derive the Baker--Campbell--Hausdorff
(BCH) formula for the product of two matrix exponentials. This formula
states that $e^{A_1} e^{A_2} = e^B$ with
\begin{equation}
\label{bch}
B = A_1 + A_2 + \tfrac12 [A_1,A_2] + \tfrac1{12} \big[A_1,[A_1,A_2]\big] +
\tfrac1{12} \big[A_2,[A_2,A_1]\big] + \cdots.
\end{equation}
Indeed, if we define the function~$A(\,\cdot\,)$ by $A(t) = A_2$ for
$t\in(0,1)$ and $A(t) = A_1$ for $t\in(1,2)$, then the solution
of~\eqref{ode} at time $t=2$ is $e^{A_1} e^{A_2}$, and the BCH
formula~\eqref{bch} is the Magnus series~\eqref{MagnExp} for this
particular choice of~$A(\,\cdot\,)$.

Magnus derived the expansion~\eqref{MagnExp} in the context of quantum
mechanics, where $A(t)$ is a skew-Hermitian matrix. Because the Magnus
series is contructed from commutators, $\Omega(t)$~is also
skew-Hermitian and $e^{\Omega(t)}$ is unitary, just like the
fundamental solution of the original differential
equation~\eqref{ode}. More generally, if $A(t)$ is in some Lie
algebra, then $\Omega(t)$ will be in the same Lie algebra and
$e^{\Omega(t)}$ will be in the corresponding Lie group.

In the 1990s, Arieh Iserles and Syvert N\o{}rsett were among the group
of mathematicians that established the discipline of \emph{Geometric
Integration}. This field concerns itself with numerical integrators
that respect the geometric structure of differential equations (see,
e.g., Hairer, Lubich and Wanner~\cite{hairer02gni} for an
introduction). Iserles and N\o{}rsett were looking for a way to
integrate~\eqref{ode} such that the numerical solution evolves on the
Lie group if the matrix~$A(t)$ is in the Lie algebra. Unaware of
Magnus' work, they rederived the Magnus series. They realized that
after truncating the infinite series and approximating the integrals
by quadrature, a method arises that respects the Lie-algebraic
structure of the equation (see Iserles and
N\o{}rsett~\cite{iserles99oso} and Iserles~\cite{iserles02etg} for
details).

Since the Magnus expansion is an infinite series, it is natural to ask
whether it converges. Indeed, Magnus himself gave an example in which
the series does not converge~\cite{magnus54oes}. The question of
convergence is the subject of this paper. The main result
(Theorem~\ref{th:conv}) gives a sufficient condition for convergence.

\section{Previous results on the convergence}
\label{sec:prev}

Due to the complexity of the expansion~\eqref{MagnExp}, several proof
strategies have been used to derive bounds on the terms $\Omega_n$,
which in turn have led to many different convergence estimates. Magnus
\cite{magnus54oes} gave no convergence estimate but stated that for
sufficiently small~$t$ the series converges.  By the known bound
$$
\|Y(t)-\Id\|_2 \leq \exp \left( \int_0^t \|A(\tau)\|_2 \,d\tau \right) - 1
$$ 
one easily arrives at the conclusion that the Magnus series converges 
whenever $\int_0^t \|A(\tau)\|_2 \,d\tau < \log 2$.%
\footnote{V.S. Varadarajan~\cite[p.~119]{varadarajan74lgl} gives
  implicitly the same result for the BCH formula.} 

However there are several improvements on this bound.  In the
following, $r$ denotes a number for which the following statement
holds:
$$
\text{If $\displaystyle \int_0^t \|A(\tau)\|_2 \, d\tau < r$ then the
Magnus expansion converges.}
$$

In the field of quantum physics there has been some interest in the
convergence issue.  In 1966, Pechukas and Light~\cite{pechukas66oef}
consider particular quantum systems and find convergence conditions,
althought these are not of the general form we are considering here
(see also Fern\'andez~\cite{fernandez90com}, Klarsfeld and
Oteo~\cite{klarsfeld89bfa}, and Salzman~\cite{salzman87ncf}).

In 1976, Karas\"ev and Mosolova~\cite{karasev77ipa} cite a bound
$r=\frac12\log2$.%
\footnote{In 1977, Suzuki~\cite{suzuki77oco} derives this bound for BCH.}
Agra{\v c}hev and Gamkrelidze, working in the field of control theory,
mention in their 1981 paper~\cite{agrachev81caa} a result by
Vakhrameev stating that $r=1.08688$.%
\footnote{Newman et al.~\cite{newman89cdf} establish the same bound
for the BCH formula in 1988.}  In 1998, this bound was rediscovered
independently by Blanes, Casas, Oteo and Ros~\cite{blanes98maf} and by
Moan~\cite{moan98eao}, using different methods. Vela~\cite{vela03aac}
states in 2003 that this is a sharp result.

In 1991, Chacon and Fomenko~\cite{chacon91rff} found an alternative
expression for~$\Omega_n$. They used this expression to prove that
$r=0.57745$.

A few years earlier, in 1987, Strichartz~\cite{strichartz87cfa} had
rediscovered the explicit expression for $\Omega_n$ found by
Bialynicki-Birula et al.~\cite{bialynicki-birula69eso}, stated in
terms of Lie brackets. He used this to prove $r=1$. The same result
was found independently by Vinokurov~\cite{vinokurov97eso} in
1997. Finally, Moan and Oteo~\cite{moan01coe} derived the bound $r=2$
(the best result at the moment) by similar techniques, except that
they avoided the use of commutators as they seemed to introduce
unneccesary complications in the convergence bound.

Moan~\cite{moan02obe} found a condition for existence of a real
logarithm for real $A$ with $r=\pi$. It was however unclear if this
condition is sufficient for convergence of the series
expansion. Theorem~\ref{th:conv} answers this question affirmatively.

\section{The existence of a real logarithm}

The fundamental solution $Y(t)$ of the differential
equation~\eqref{ode} is an invertible matrix. A theorem by
Gantmacher~\cite{gantmacher59tom} states that every invertible matrix
has a logarithm, and hence there exists an~$\Omega(t)$ such that
$e^{\Omega(t)} = Y(t)$. However, the logarithm may fail to be real,
even if $Y(t)$ is real. For example, the matrix
\begin{equation}
\left[\begin{array}{cc} -1 & 0 \\ b & -1 \end{array}\right],
\qquad b \in \RR\setminus\{0\},
\label{nolog} 
\end{equation}
does not have a real logarithm. Since all the terms in the Magnus
expansion~\eqref{MagnExp} are real if the original differential
equation~\eqref{ode} is real, we conclude that the infinite series
cannot converge to a logarithm of~$Y(t)$ if $Y(t)$ has no real
logarithm. The question we are studying in this section is therefore:
Does $Y(t)$ have a real logarithm?

The following lemma is easily proved, for instance by factoring
$\Phi=V J V^{-1}$ where $J$ is in Jordan form.

\begin{lemma}\label{lemmalog} 
Suppose that the invertible matrix $\Phi$ has no negative eigenvalues,
that is, suppose that the eigenvalues of~$\Phi$ are contained in
$\CC \setminus (-\infty,0]$. Then
\begin{equation}
\label{matrixlog}
\log \Phi = (\Phi-\Id) \int_0^\infty \frac{1}{1+\mu}
(\mu\Id+\Phi)^{-1} \, d\mu.
\end{equation}
\end{lemma}

\noindent
It follows that if $Y(t)$ is real and has no negative eigenvalues,
then the logarithm of~$Y(t)$ is real as well. The next result (which
can be found in~\cite{moan02obe}) gives an easy condition under which
$Y(t)$ has no negative eigenvalues and hence a real logarithm.

\begin{theorem} \label{theoexist}
Let $A(t)$ be a real integrable matrix, and let $Y(t)$ denote the
solution of $Y' = A(t)\,Y$, $Y(0)=\Id$. If $\int_0^t \|A(\tau))\|_2
\,d\tau < \pi$, then $Y(t)$ has a real logarithm.
\end{theorem}

\begin{proof} 
Choose an arbitrary vector~$y_0$, and consider the vector $y(t)=
Y(t)\,y_0$ satisfying $y' = A(t)\,y$ and $y(0)=y_0$. Let $\hat{y}(t) =
y(t) / \|y(t)\|_2$ denote the unit vector in the direction~$y(t)$.
Then 
$$
y'(t) = \big( \tfrac{d}{dt} \|y(t)\|_2 \big) \hat{y}(t) + 
\|y(t)\|_2 \, \hat{y}'(t).
$$
Since $\hat{y}(t)$ and $\hat{y}'(t)$ are orthogonal, it follows that
$$
\|y'(t)\|_2^2 = \big( \tfrac{d}{dt} \|y(t)\|_2 \big)^2 +
\|\hat{y}'(t)\|_2^2 \, \|y(t)\|_2^2.
$$
Therefore, we have
$$
\|\hat{y}'(t)\|_2 \, \|y(t)\|_2  \le \|y'(t)\|_2 \le \|A(t)\|_2 \|y(t)\|_2,
$$
and thus $\|\hat{y}'(t)\|_2 \le \|A(t)\|_2$. Integrating this
inequality, we find the bound  
$$
\int_0^t \|\hat{y}'(\tau)\|_2 \,d\tau \le \int_0^t \|A(\tau)\|_2 \,d\tau.
$$
The left-hand side is the length of the curve swept out by the unit
vector $\hat{y}'(\tau)$ when $\tau \in [0,t]$. Therefore, the angle
between~$y(t)$ and~$y(0)$ is smaller than~$\pi$ if $\int_0^t
\|A(\tau))\|_2 \,d\tau < \pi$. Since $y(t) = Y(t)\,y(0)$, this implies
that $Y(t)$ has no negative eigenvalues, and hence it has a real
logarithm.
\end{proof}

\noindent
Example~\ref{ex:pcthesis} in the forthcoming
Section~\ref{sec:examples} shows that the constant~$\pi$ in
Theorem~\ref{theoexist} is sharp.

\section{Proof of convergence}
\label{sec:proof}

Theorem~\ref{theoexist} gives a condition for $Y(t)$ to have a real
logarithm. The following theorem states that under the same condition,
the Magnus series converges to this logarithm.

\begin{theorem}
\label{th:conv}
Let $A(t)$ be a real integrable matrix, and let $Y(t)$ denote the
solution of $Y' = A(t)\,Y$, $Y(0)=\Id$. If $\int_0^t \|A(\tau)\|_2
\,d\tau < \pi$, then the Magnus series~\eqref{MagnExp} converges and
its sum $\Omega(t)$ satisfies $e^{\Omega(t)} = Y(t)$.
\end{theorem}

\begin{proof}
We write the Magnus series as $\Omega(t) = \sum_{n=1}^\infty
\Omega_n(t)$, where the term $\Omega_n(t)$ is a sum of $n$-fold
integrals of $n-1$ nested commutators. If we now introduce a new
parameter~$\kappa$ and we replace $A(t)$ by~$\kappa A(t)$, then
the Magnus series becomes $\Omega(t;\kappa) = \sum_{n=1}^\infty
\kappa^n \Omega_n(t)$. The idea is to fix~$t$ and consider the
function~$f$ defined by $f(\kappa) = \log Y(t;\kappa)$ where
$Y(t;\kappa)$ denotes the solution of $Y' = \kappa A(t) Y$, $Y(0) =
\Id$ with $\kappa\in\CC$. We will show that $f$ is analytic and that
the Magnus series is the Taylor series of this function around
$\kappa=0$. We can then use the standard result from the theory of
complex functions which states that the Taylor series of a function
converges in a disc if the function is analytic in that disc.

Set $\gamma = \int_0^t \|A(\tau)\|_2 \,d\tau$. As stated in
Section~\ref{sec:prev}, it is easy to show that if $|\kappa| <
\frac1\gamma \log 2$, the Magnus series converges and its sum
$\Omega(t;\kappa)$ satisfies $e^{\Omega(t;\kappa)} =
Y(t;\kappa)$. Hence, the power series $\Omega(t;\kappa)$ coincides
with $f(\kappa)$ for $|\kappa| < \frac1\gamma \log 2$, and the Magnus
series is the Taylor series expansion of~$f$ around $\kappa = 0$.

We say that a matrix-valued function (like~$f$) is analytic if all the
matrix entries are analytic functions. We now want to prove that $f$
is analytic in the disc with radius~$\frac\pi\gamma$. Firstly, the
fundamental matrix~$Y(t;\kappa)$ is analytic as a function
of~$\kappa$. The proof of Theorem~\ref{theoexist} shows that
$Y(t;\kappa)$ has no eigenvalues in $(-\infty,0]$ if $|\kappa| <
\frac\pi\gamma$, so the logarithm is given by~\eqref{matrixlog}.
Hence, the derivative of $f$ is given by
\begin{multline*}
f'(\kappa) = \frac{\partial}{\partial\kappa} Y(t;\kappa) 
\int_0^\infty \frac1{\mu+1} \big(\mu\Id+Y(t;\kappa)\big)^{-1} \,d\mu
\\
{} - Y(t;\kappa) \int_0^\infty \frac1{\mu+1}
\big(\mu\Id+Y(t;\kappa)\big)^{-1} \frac{\partial}{\partial\kappa}
Y(t;\kappa) \big(\mu\Id+Y(t;\kappa)\big)^{-1} \,d\mu.
\end{multline*}
The right-hand side is well-defined, thus proving that the
function~$f$ is analytic in the disc with radius~$\frac\pi\gamma$.
Hence, we can expand the entries of the matrix~$f(\kappa)$ in a power
series, and this series will converge provided that $|\kappa| <
\frac\pi\gamma$. But this power series is precisely the Magnus series.
\end{proof}

\section{Examples}
\label{sec:examples}

In this section, we study some examples to investigate the connections
between the condition $\int_0^t \|A(\tau)\|_2 \,d\tau < \pi$, the
eigenvalues of the fundamental solution, the existence of a real
logarithm, and the convergence of the Magnus series.

\begin{example}
\label{ex:constant}
The following simple example suffices to show that the condition
$\int_0^t \|A(\tau)\|_2 \,d\tau < \pi$ is not a necessary
condition. Consider the equation $Y' = AY$ where $A = \bigl[
\begin{smallmatrix} 0 & 1 \\ -1 & 0 \end{smallmatrix} \bigr]$ is a
constant matrix. The fundamental matrix is
$$
Y(t) = \begin{bmatrix}
\cos t & -\sin t \\ \sin t & \cos t
\end{bmatrix}.
$$ 
When $t=\pi$, this matrix has a double eigenvalue at~$-1$.
Furthermore, we have $\int_0^\pi \|A(t)\|_2 \,dt = \pi$.  Nevertheless,
$Y(t)$~has a real logarithm for all~$t$, including $t=\pi$, because
$Y(t) = e^{tA}$. The Magnus series converges for all~$t$: the terms in
the series are $\Omega_1=A$ and $\Omega_k=0$ for $k>1$. The critical
point here seems to be that the double eigenvalue at~$-1$ is not
defective (meaning that its algebraic multiplicity equals its
geometric multiplicity).
\end{example}

\begin{example}
\label{ex:pcthesis}
This example, taken from Moan~\cite{moan02obe}, shows that the condition
$\int_0^t \|A(\tau)\|_2 \,d\tau < \pi$ is sharp, in the sense that the
constant~$\pi$ on the right-hand side cannot be replaced by a bigger
constant. Consider the equation $Y' = A(t)Y$ with
$$
A(t) = \frac12 \begin{bmatrix}
\sin 2t & -1-\cos 2t \\ 1-\cos 2t & -\sin 2t 
\end{bmatrix}.
$$
A simple computation shows that $\int_0^\pi \|A(\tau)\|_2 \,d\tau =
\pi$, and that the solution of the differential equation is given by
$$
Y(t) = \begin{bmatrix}
t \sin t + \cos t & -\sin t \\ 
\sin t -t \cos t & \cos t
\end{bmatrix}.
$$
Hence $Y(\pi) = \big[ \begin{smallmatrix} -1 & 0 \\ \pi & -1
\end{smallmatrix} \big]$, which is of the form~\eqref{nolog}, and
therefore $Y(\pi)$ does not have a real logarithm. This implies
that the Magnus series diverges for $t=\pi$.
\end{example}

\begin{example}
\label{ex:powerseries}
This example shows that the Magnus series may diverge even though the fundamental
matrix~$Y(t)$ has a real logarithm. Take 
$$ 
A(t) = \begin{bmatrix} 2 & t \\ 0 & -1 \end{bmatrix}.
$$
The Magnus series is
$$
\Omega(t) = \begin{bmatrix} 2t & \frac12t^2 \\ 0 & -t \end{bmatrix} 
+ \begin{bmatrix} 0 & -\frac14t^3 \\ 0 & 0 \end{bmatrix}
+ \begin{bmatrix} 0 & 0 \\ 0 & 0 \end{bmatrix}
+ \begin{bmatrix} 0 & \frac1{80}t^5 \\ 0 & 0 \end{bmatrix} + \cdots
$$
On the other hand, the solution of $Y' = A(t) Y$ is
$$
Y(t) = \begin{bmatrix}
e^{2t} & \frac19 e^{2t} - \bigl( \frac19 + \frac13 t \bigr) e^{-t} \\
0 & e^{-t}
\end{bmatrix}.
$$
The logarithm of this solution is
$$
\log Y(t) = \begin{bmatrix}
2t & f(t) \\ 0 & -t
\end{bmatrix} \quad\text{where}\quad 
f(t) = \frac{te^{2t}-(t+3t^2)e^{-t}}{3(e^{-2t}-e^t)} \, .
$$
So, the Magnus series is the Taylor series expansion of $\log Y(t)$
around $t=0$ (this is not true in general). The function~$f$ has a
pole at $t = \frac23\pi i$, thus the Magnus series converges up to $t
= \tfrac23\pi$. 

However, $Y(\tfrac23\pi)$ has eigenvalues at~$e^{2\pi/3}$
and~$e^{4\pi/3}$, so $Y(\tfrac23\pi)$ has a real logarithm.
Nevertheless, the Magnus series diverges. The reason is probably as
follows. As in the proof of Theorem~\ref{th:conv}, let $Y(t;\kappa)$
denote the fundamental solution of $Y' = \kappa A(t) Y$. The
eigenvalues of~$Y(t;\kappa)$ are $e^{\kappa t}$ and~$e^{-2\kappa
t}$. When $\kappa = i$, these eigenvalues move in a circle around the
origin and collide when $t = \tfrac23\pi$. This collision causes the
Magnus series to diverge at $t = \tfrac23\pi$.
\end{example}

\begin{example}
\label{ex:compl}
The preceding examples all have some structure which allowed us to
determine where the Magnus series starts to diverge. We close this
section with a more-or-less randomly chosen example.

Consider the equation $Y' = A(t) Y$, where the matrix~$A$ is defined by
\begin{equation}
\label{complex}
A(t) = \begin{bmatrix}
-t & 3t & 0 & -3t^2+t+2 \\ t^2-t & -3 & t^2+2t+3 & 0 \\
3 & 0 & t^2-2t & -t^2-3 \\ t^2-t+3 & 2t^2-3t & -3t-2 & -t+2
\end{bmatrix}.
\end{equation}
Since $A$ is polynomial in~$t$, all the commutators in the Magnus
series are also polynomials and hence the integrals can be computed
exactly. The first two terms of the series are given by
\begin{align*}
\Omega_1 &= \begin{bmatrix}
-\frac12 t^2 & \frac32 t^2 & 0 & -t^3 + \frac12 t^2 + 2t \\
\frac13 t^3 - \frac12 t^2 & -3t & \frac13 t^3 + t^2 + 3t & 0 \\
3t & 0 & \frac13 t^3 - t^2 & -\frac13 t^3 - 3t \\
\frac13 t^3 - \frac12 t^2 + 3t & \frac23 t^3 - \frac32 t^2 
    & -\frac32 t^2 - 2t & -\frac12 t^2 + 2t
\end{bmatrix} \\
\Omega_2 &= \begin{bmatrix}
  \scriptstyle -\frac1{60} t^5 - \frac{11}{12} t^4 + \frac5{12} t^3 
& \scriptstyle \frac7{60} t^5 - \frac13 t^4 - \frac14 t^3
& \scriptstyle \frac1{10} t^5 + \frac12 t^4 + \frac{13}{12} t^3
& \scriptstyle -\frac12 t^4 + \frac16 t^3 \\
  \scriptstyle -\frac1{60} t^5 + \frac12 t^4 + \frac14 t^3
& \scriptstyle \frac1{20} t^5
& \scriptstyle -\frac1{15} t^5 + t^3
& \scriptstyle \frac16 t^4 - \frac23 t^3 \\
  \scriptstyle \frac1{60} t^5 + \frac14 t^4 - \frac12 t^3 
& \scriptstyle \frac1{20} t^5 + \frac12 t^4 - \frac32 t^3
& \scriptstyle \frac1{20} t^5 + \frac16 t^4 - \frac34 t^3
& \scriptstyle -\frac1{60} t^5 + \frac13 t^4 \\
  \scriptstyle \frac1{60} t^5 - \frac16 t^4 - \frac7{12} t^3
& \scriptstyle \frac1{12} t^5 - \frac56 t^4 + \frac12 t^3
& \scriptstyle \frac16 t^5 + \frac23 t^4 - \frac5{12} t^3
& \scriptstyle -\frac1{12} t^5 + \frac34 t^4 + \frac13 t^3
\end{bmatrix}
\end{align*}
We computed the first thirty terms of the Magnus series with the help
of a computer algebra system. The recursive formulas given by Blanes,
Casas, Oteo and Ros~\cite{blanes98maf} proved to be useful for this
purpose; other formulations require a very long time to evaluate. In
Figure~\ref{fig1}, we plot the sums of the first fifteen, twenty,
twenty-five and thirty terms of the Magnus series~\eqref{MagnExp}. It
seems that the series starts to diverge between $t = 0.7$ and $t =
0.8$, though it is of course not possible to pinpoint the location
precisely.

\begin{figure}
  \begin{center}
    \includegraphics[width=1.0\linewidth]{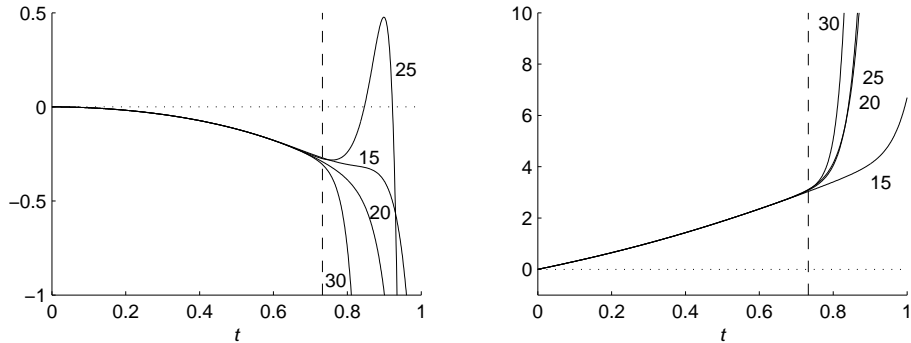}
  \end{center}
  \caption{The left plot shows the partial sums of the first 15, 20,
  25 and 30 terms of the (1,1)~element in the Magnus
  expansion~\eqref{MagnExp}, when $A$ is as given in~\eqref{complex}.
  The dash line indicates $t = 0.733$. The right panel shows the same
  for the (2,3)~element.}
  \label{fig1}
\end{figure}

Figure~\ref{fig2} shows how the eigenvalues of the fundamental
solution~$Y(t)$ move around as $t$ increases from~0 to~1. The
eigenvalues do not collide in the left half-plane, so $Y(t)$ has a
real logarithm for $t\in[0,1]$. However, as the right panel in the
figure shows, the eigenvalues of $Y(t;e^{i\alpha_*})$ with $\alpha_* =
1.805\ldots$ do collide. The collision takes place at $t_* =
0.733\ldots$ and $\lambda_* \approx -0.485 + 0.0249i$. The value
of~$t_*$ is shown by the dash line in Figure~\ref{fig1}. It seems
plausible that the Magnus series starts to diverge around this point.

\begin{figure}
  \begin{center}
    \includegraphics[width=1.0\linewidth]{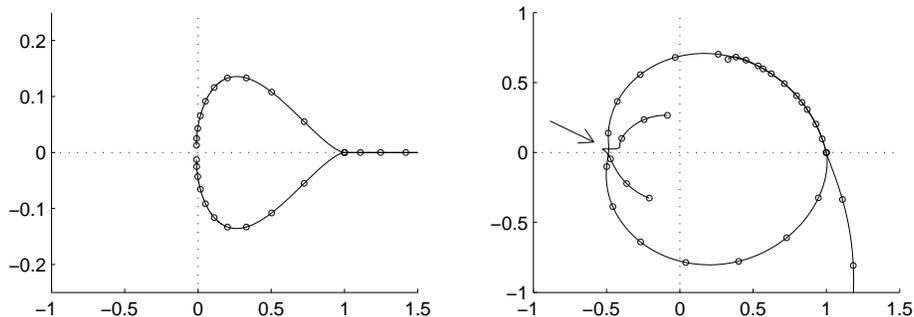}
  \end{center}
  \caption{The left plot shows the location in the complex plane of
  the eigenvalues of the fundamental solution of $Y' = A(t) Y$ with
  $A$ as in~\eqref{complex}, for $t \in [0,1]$. All four eigenvalues
  are at 1 when $t=0$. As $t$ increases, two eigenvalues move to the
  right, while the other two form a complex pair. The small circles
  show the locations at $t = 0, 0.1, 0.2, \ldots, 1$. It is clear that
  no eigenvalues collide on the negative real axis for $t \in
  [0,1]$. \protect\\ 
  The right plot shows the eigenvalues for $Y' = e^{i\alpha_*}
  A(t) Y$ with $\alpha_* \approx 1.805$. Two eigenvalues encircle the
  origin and collide near the arrow.}  
  \label{fig2}
\end{figure}

Incidentally, the condition $\int_0^t \|A(\tau)\|_2 \,d\tau < \pi$ is
satisfied for $t < 0.56\ldots$. However, the Magnus series continues
to converge for slightly larger values of~$t$, showing again that this
condition is not necessary for convergence. At $t=t_*$, the integral
is approximately~$4.36$.
\end{example}

\section{Conclusion}

The main result of the paper is Theorem~\ref{th:conv}, which states
that the Magnus series converges if $\int_0^t \|A(\tau)\|_2 \,d\tau <
\pi$. This condition is in the same form as various earlier results
mentioned in Section~\ref{sec:prev}. Example~\ref{ex:pcthesis} shows
that our result is sharp in the sense that the constant~$\pi$ is the
largest number for which the result holds. However, the other examples
show that the condition is not necessary for convergence.

The proof of Theorem~\ref{th:conv} shows that the Magnus series can be
considered as the Taylor series of $\log Y(t;\kappa)$ around $\kappa =
0$, and hence the radius of convergence is determined by the nearest
singularity. Lemma~\ref{lemmalog} implies that there are no
singularities if the eigenvalues of $Y(t;\kappa)$, which start at 1
when $\kappa = 0$, do not cross the negative real axis as $\kappa$
moves in the unit disc. However, the choice of the negative real axis
is arbitary; a formula similar to~\eqref{matrixlog} holds for any
other branch cut of the logarithm.

Examples~\ref{ex:powerseries} and~\ref{ex:compl} suggest that
divergence of the Magnus series is associated with eigenvalues
of~$Y(t;\kappa)$ encircling the origin and colliding as $\kappa$ moves
from~0. In this situation, the eigenvalues are on different sheets of
the Riemann surface of the logarithm when they collide. The authors
feel that this conflict leads to the divergence of the Magnus
series. However, as example~\ref{ex:constant} shows, not every
collision of eigenvalues leads to divergence of the Magnus series. If
the multiple eigenvalue at the collision is not defective, then the
eigenvalues retain their identity throughout the process and no
conflict arises.

Of course, this is not a proof, but we think that it may be possible
to make it rigorous using the correct formalism. We are thus led to
propose the following conjecture, which gives a necessary and
sufficient condition for convergence.

\begin{conjecture}
Let $Y(t;\kappa)$ denote the solution of $Y' = \kappa A(t) Y$,
$Y(0)=I$. Denote the eigenvalues of $Y(t;\kappa)$ by
$\lambda_n(t;\kappa)$, where the eigenvalues are to be ordered so that
$\lambda_n$ is a continuous function of~$t$. Let $t_*$ be the smallest
$t>0$ for which there exists a $\kappa\in\CC$ with $|\kappa| = 1$ such
that there is a multiple eigenvalue, say $\lambda_i(t,\kappa) =
\lambda_j(t,\kappa)$ with $i \ne j$, for which the geometric
multiplicity is smaller than the algebraic multiplicity, and the loop
$$
\{ \lambda_i(\tau,\kappa) \mid \tau\in[0,t] \} \cup
\{ \lambda_j(\tau,\kappa) \mid \tau\in[0,t] \}
$$
encircles the origin. Then the Magnus series converges if and only if
$t<t_*$.
\end{conjecture}

\noindent
This conjecture suggests two tasks. The first is obviously whether we
can find a proof for this conjecture. However, even if the conjecture
is true, it is difficult to apply in practice because the condition is
not easy to check. Therefore, the second task is to find a more
practical condition for convergence.

\medskip\noindent
\textbf{Acknowledgements.} We thank Arieh Iserles for introducing us
to the Magnus series, and for being a pleasant and enthusiastic mentor
for us. JN was supported by EPSRC First Grant GR/S22134/01.

\bibliographystyle{abbrv}
\bibliography{../jitse}

\end{document}